\magnification 1200
\font\bb=msbm10

\def\N{{\hbox{\bb N}}}
\def\norm#1{{\|#1\|}}
\def\Spec#1{\Phi _ {#1}}
\def\qed{{\hfill\vrule height 4pt width 4pt depth 0pt
             \par\vskip\baselineskip}}

\centerline{\bf STRONG DITKIN ALGEBRAS WITHOUT}
\smallskip
\centerline{\bf BOUNDED RELATIVE UNITS}
\medskip
\centerline{J.F.~FEINSTEIN}
\medskip
\centerline{School of Mathematical Sciences, Pure Mathematics Division,}
\smallskip
\centerline{University of Nottingham, Nottingham NG7 2RD,
England}
\smallskip
\centerline{e-mail: Joel.Feinstein@nottingham.ac.uk}
\medskip
\noindent
{\bf Abstract.} In a previous note the author  
gave an example of a strong Ditkin algebra which does not have
bounded relative units in the sense of Dales. In this note we 
investigate a certain family of Banach function algebras on the 
one point compactification of \N, and see that within this 
family are many easier examples of strong Ditkin algebras without 
bounded relative units in the sense of Dales.
\medskip
\centerline{1. INTRODUCTION}
\medskip
Regularity conditions for Banach function algebras have important applications 
in several areas of functional analysis, including automatic continuity theory
and the theory of Wedderburn decompositions (see, for example, [{\bf 3}]). 
There is also a close 
connection between regularity and the theory of decomposable operators,
as was shown by Neumann in [{\bf 6}]. It is thus important both to investigate 
the connections between these regularity conditions and other conditions that
a Banach function algebra may satisfy, and also to find out what the 
relationships are  between these regularity conditions. 
For a survey of many
regularity conditions and their relationships to each other, see [{\bf 5}].
\smallskip
The aim of this note is to investigate some of the stronger regularity
conditions for a particular class of Banach function algebras, following on
from the work in [{\bf 4}].
In [{\bf 4}], the author solved a problem posed by Bade (see [{\bf 2}]), 
by showing
that every strong Ditkin algebra has bounded relative units at each point
of its character space. The proof was very short and elementary. In fact,
unknown to the author, essentially the same argument had been used earlier
by Bachelis, Parker and Ross ([{\bf 1}]) to prove a special case of the same
result. (I am grateful to Jan Stegeman for later pointing out the
existence of this paper to me). Also in [{\bf 4}], the author constructed an
example of a strong Ditkin algebra which does not satisfy the stronger
condition of having bounded relative units in the sense of Dales (see
below). In Section 2 we shall introduce a particular class of Banach
function algebras on $\N_\infty$, and determine precisely which combinations
of regularity properties are possible for algebras in this class.
In particular we shall show that this class includes many further
examples of strong Ditkin algebras which do not have
bounded relative units in the sense of Dales. These examples are easier
than the example constructed in [{\bf 4}], and show that interesting 
combinations of properties are possible even when working with Banach
function algebras on $\N_\infty$.
First we recall the 
standard notation and definitions which we shall need.
\smallskip
\noindent
{NOTATION.} For any compact, Hausdorff space $X$ we denote the 
algebra of all continuous, complex-valued functions on $X$ by $C(X)$. 
Let $A$ be a commutative,
unital algebra. We denote by $\Spec{A}$ the character space of $A$.
Let $S$, $T$ be subsets of $A$. We denote by $S \cdot T$ the set
$$\{ab:a \in S, b \in T \}.$$
\par
We denote the one point compactification of $\N$ by $\N_\infty$ (so that
$\N_\infty = \N \cup \{\infty\}$).
For $f\in C(\N_\infty)$, $\norm{f}_\infty$ denotes the usual uniform norm of $f$ 
on $\N_\infty$, so that
$$\norm{f}_\infty= \sup \{|f(n)|: n\in\N\}.$$
\medskip
{DEFINITION 1.1:} Let $X$ be a compact space. A {\it function algebra} on $X$
is a subalgebra of $C(X)$ which contains all of the constant functions and 
which separates the points of $X$. 
A {\it Banach function algebra} on $X$ is a function algebra on $X$ with a 
complete algebra norm.
\smallskip
Throughout we shall consider only unital Banach algebras. Using the Gelfand
transform, every commutative, semisimple Banach algebra may be regarded
as a Banach function algebra on its character space.
\smallskip
\noindent
{NOTATION.} Let $A$ be a (unital) Banach function algebra on $\Spec{A}$
and let $E$ be a closed subset of $\Spec{A}$. We define two ideals $J(E)$,
$I(E)$ in $A$ as follows:
$$J(E)=\{f \in A: f \hbox{ vanishes on some neighbourhood of } E\};$$ 
$$I(E)=\{f \in A: f(E) \subseteq \{0\}\}.$$ 
For $\phi \in \Spec{A}$, we set $M_\phi = I(\{\phi\})$ and
we denote $J(\{\phi\})$ by $J_\phi$.
\medskip
{DEFINITION 1.2:} Let $A$ be a Banach function algebra on $\Spec{A}$ and let
$\phi \in \Spec{A}$. Then $A$ is {\it strongly regular at} $\phi$ if
$J_\phi$ is dense in $M_\phi$; 
$A$ satisfies {\it Ditkin's condition 
at\/} $\phi$ if $J_\phi \cdot M_\phi$ is dense in $M_\phi$.
We say that $A$ has {\it bounded relative units at \/} 
$\phi$ if there exists $C \geq 1$ such that, for all compact sets
$E \subseteq \Spec{A} \backslash \{ \phi \} $, there is $f \in J_\phi $ with 
$\norm{f} \leq C$ such that $f(E) \subseteq \{ 1 \}.$
The algebra $A$ is {\it strongly regular\/} if it is strongly regular at every 
point of $\Spec{A}$; $A$ is a {\it Ditkin algebra\/} if it satisfies Ditkin's 
condition at 
every point of $\Spec{A}$; $A$ is a {\it strong Ditkin algebra\/} if 
every maximal ideal of $A$ has a bounded approximate identity and $A$ is 
strongly regular. The Banach function algebra $A$ {\it has spectral
synthesis\/} if, for every closed set $E \subseteq \Spec{A}$, $J(E)$ is
dense in $M(E)$.
\smallskip
Clearly every Ditkin algebra is strongly regular. In the special case where 
$\Spec{A} = \N_\infty$ it is easy to see that if $A$ is a Ditkin algebra, then
$A$ has spectral synthesis (this is not true in general).
\smallskip
There are two definitions for a Banach function algebra $A$ to have 
bounded relative units. One is that, for every $\phi \in \Spec{A}$, $A$ has
bounded relative units at $\phi$. If this condition holds we shall say that
$A$ has {\it bounded relative units in the sense of Bade}. The other 
definition, used more frequently, is stronger, insisting
that the constant $C$ involved does not depend on $\phi$. If this stronger
condition holds we shall say that $A$ has {\it bounded relative units in
the sense of Dales}. For uniform algebras, these two conditions are 
equivalent. But they differ for general Banach function algebras, as was 
shown in [{\bf 4}]. 

The following result combines some standard theory with [{\bf 4}, Theorem 5].

{\bf PROPOSITION 1.3.} {\sl Let $A$ be a Banach function algebra. Then $A$ is a
strong Ditkin algebra if and only if $A$ is strongly regular and has
bounded relative units in the sense of Bade.}
\smallskip
As we mentioned earlier,  in [{\bf 4}] there is an example of a
strong Ditkin algebra which does not have
bounded relative units in the sense of Dales. In the next section we shall
give some easier examples of this.
\medskip
\centerline{2. BANACH FUNCTION ALGEBRAS ON $\N_\infty$.}
\medskip
We shall now introduce the class of Banach function algebras that we shall work
with. From now on we shall use $\alpha$ to denote a sequence of positive real 
numbers, with $\alpha = (\alpha_n)_{n=1}^\infty$. Given such a sequence 
$\alpha$, we define $A_\alpha$ by
$$A_\alpha=\left\{{f\in C(\N_\infty):\sum_{n=1}^\infty
{\alpha_n|f(n+1)-f(n)|}<\infty}\right\}.$$
It is easy to see that $A_\alpha$ is a subalgebra of $C(\N_\infty)$, and
that $A_\alpha$ is a Banach function algebra, where the norm of a function
$f \in A_\alpha$ is given by
$$\norm{f}= \norm{f}_\infty + \sum_{n=1}^\infty{\alpha_n|f(n+1)-f(n)|}.$$
It is also easy to see that the character space of $A_\alpha$ is
just $\N_\infty$ (this follows from [{\bf 4},~Proposition~6], 
for example).\par
We shall now investigate the further properties of these algebras.
We shall see that $A_\alpha$ is always a Ditkin algebra (and hence
strongly regular). In fact we will see that $A_\alpha$ has
a stronger property: for all $x \in \N_\infty$, $M_x$ has an
approximate identity (not necessarily bounded) consisting of a
sequence of elements of $J_x$. Of course, for each $x\in \N$,
it is clear that $J_x = M_x$, and that these ideals have
an identity. So we only need to check what happens at the point
$\infty$. We begin by looking at the particular sequence of functions 
$(e_k) \subseteq A_\alpha$ 
defined by
$$e_k(n)=\cases{1&if $n\leq k$,\cr
0&otherwise.}$$
We shall see that this sequence of functions in $J_\infty$ 
always has a subsequence which is an approximate identity for 
$M_\infty$, although the sequence $(e_k)$ need not itself be
such an approximate identity. First we need an elementary lemma.
\smallskip
{\bf LEMMA 2.1.} {\sl Let $(n_k)_{k=1}^\infty$ be a strictly increasing
sequence of natural numbers. Then the following conditions are equivalent:
\item{(a)} $(e_{n_k})$ is an approximate identity for $M_\infty$;
\item{(b)} for all $f \in M_\infty$, $\lim_{k\to\infty} \alpha_{n_k} f(n_k + 1) = 0$;
\item{(c)} for all $f \in M_\infty$, $\lim_{k\to\infty} \alpha_{n_k} f(n_k) = 0$.
}
\par
{PROOF:} For any $f \in A_\alpha$ we know that 
$\sum_{k=1}^\infty{\alpha_n|f(n+1)-f(n)|}<\infty$, and so
$\lim_{n \to \infty} \alpha_n|f(n+1)-f(n)| = 0$. From this it is clear that
(b) and (c) are equivalent. \par
Now let $f \in M_\infty$.
Then $$(f - e_{n_k} f)(j)=\cases{0&if $j\leq n_k$,\cr
f(j)&otherwise.}$$
We have that
$$\norm{f - e_{n_k} f} = \norm{f - e_{n_k} f}_\infty + 
\sum_{j=n_k + 1}^\infty{(\alpha_j|f(j+1)-f(j)|)} + \alpha_{n_k}|f(n_k + 1)|.\eqno(2.1)$$
Since the first two terms on the right hand side 
of (2.1) tend to zero as $k \to \infty$, it is clear that
$\lim_{k \to \infty} \norm {f - e_{n_k}f} = 0$ if and only if 
$\lim_{k \to \infty}\alpha_{n_k} f(n_k + 1) = 0$. 
It is now immediate that (a) and (b) are equivalent, as required. \qed
\medskip
The sequence $(e_k)$ is not always an approximate identity for
$M_\infty$. For example, suppose that $\alpha$ is the sequence
$(\alpha_n)$ defined by
$$\alpha_n = \cases{1&if $n$ is odd,\cr
{n\over 2}&if $n$ is even.}$$
Then $(e_k)$ is not an approximate identity for $M_\infty$. This can be seen by 
considering the function $f \in C(\N_\infty)$ which satisfies
$$f(j) = 2^{-k}$$
whenever $j, k \in \N$ are such that
$2^{k-1} \leq j < 2^k$.\par
When $j = 2^k - 1$, $\alpha_j |f(j+1) - f(j)| = 2^{-k-1}$. All other
$|f(j+1)-f(j)|$ are $0$, so that $f \in A_\alpha$, and $f \in M_\infty$. But 
$\alpha_{2^k} f(2^k) = {1 \over 4}$ for all $k$, and so $\alpha_n f(n)$ does not
tend to $0$ as $n \to \infty$. Thus, by Lemma 2.1, $(e_k)$ is not an 
approximate identity for $M_\infty$. 
\medskip
We now continue our investigation into when $(e_{n_k})$ is an approximate 
identity.
\medskip
{\bf THEOREM 2.2.} {\sl Let $(n_k)$ be as in Lemma 2.1.
\item{(a)} Suppose that, for all $k \in \N$, 
$\alpha_{n_k} = \inf \{\alpha_j: j\geq n_k\}.$ Then $(e_{n_k})$ is an
approximate identity for $M_\infty$.
\item{(b)} The sequence $(e_{n_k})$ is a bounded approximate identity for
$M_\infty$ if and only if the sequence $(\alpha_{n_k})$ is bounded.}
\smallskip
{PROOF:} To prove (a), suppose that $(\alpha_{n_k})$ satisfies the 
given condition. Let $f \in M_\infty$. Then we have, for $k \in \N$,
$$f(n_k) = \sum_{j=n_k}^\infty{(f(j)-f(j+1))}.$$
Thus
$$\eqalign{|\alpha_{n_k} f(n_k)| &\leq 
\sum_{j=n_k}^\infty{\alpha_{n_k}|f(j)-f(j+1)|}\cr
&\leq \sum_{j=n_k}^\infty{\alpha_j |f(j)-f(j+1)|}.}$$
Since the last sum tends to $0$ as $k \to \infty$, it follows
that $\lim_{k\to\infty} \alpha_{n_k} f(n_k) = 0$. It now follows
from Lemma 2.1 that $(e_{n_k})$ is an approximate identity for
$M_\infty$.\par
To prove (b), first note that, for all $k$,
$$\norm{e_{n_k}} = 1 + \alpha_{n_k}.$$
Thus it is clear that the sequence $(e_{n_k})$ is bounded in norm
if and only if $(\alpha_{n_k})$ is bounded. It remains to show that,
in this case, $(e_{n_k})$ is also an approximate identity for $M_\infty$.
But this is immediate from Lemma 2.1. The result follows. \qed
\medskip
We can now see that $(e_k)$ is itself often an
approximate identity.
\smallskip
{\bf COROLLARY 2.3.} {\sl If the sequence $\alpha$ is either nondecreasing
or bounded then $(e_k)$ is an approximate identity for $M_\infty$.}
\smallskip
{PROOF:} This is immediate from Theorem 2.2. \qed
\medskip
We can now prove our main result about the existence of approximate
identities in $M_\infty$.
\smallskip
{\bf THEOREM 2.4.} {\sl The sequence $(e_k)$ always has a subsequence which is
an approximate identity in $M_\infty$ in $A_\alpha$.}
\smallskip
{PROOF:} If $(\alpha_n)$ has a bounded subsequence, then the result is immediate 
from Theorem 2.2 (b). Thus we may assume that $(\alpha_n)$ diverges to $\infty$.
But then it is easy to choose $(n_k)$ such that the conditions of Theorem 2.2. (a)
are satisfied, and so $(e_{n_k})$ is an approximate identity for $M_\infty$. \qed
\medskip
As a corollary, we obtain immediately the fact that each $A_\alpha$ is a
Ditkin algebra.
\medskip
{\bf COROLLARY 2.5.} {\sl For every sequence of positive real numbers 
$\alpha$, the Banach function algebra $A_\alpha$ is a Ditkin algebra.}
\smallskip
{PROOF:} For $x\in \N$ it is clear that $J_x = M_x$ and that these ideals have
an identity. It remains to check that $J_\infty \cdot M_\infty$ is dense in 
$M_\infty$, but this is immediate from Theorem 2.4. The result follows. \qed
\medskip
Of course this means that $A_\alpha$ is always strongly regular. By the
nature of the algebra, it follows that $A_\alpha$ is always separable.
Since $\Spec{A_\alpha} = \N_\infty$, it also follows that $A_\alpha$
has spectral synthesis.
\medskip
We now characterise those sequences $\alpha$ for which the algebra 
$A_\alpha$ has various regularity properties. 
\medskip
{\bf THEOREM 2.6.} {\sl Let $\alpha = (\alpha_n)_{n=1}^{\infty}$ be any 
sequence of positive real numbers. Then\par
\item{(a)} $A_\alpha$ is a strong Ditkin algebra if and only if 
$\liminf_{n\to\infty} \alpha_n < \infty$, i.e. $(\alpha_n)$ does not 
diverge to $\infty$.
\smallskip
\item{(b)} In $A_\alpha$, $M_\infty$ has a bounded approximate identity if and 
only if $\liminf_{n\to\infty} \alpha_n < \infty$.
\smallskip
\item{(c)} The Banach function algebra $A_\alpha$ has bounded relative units
in the sense of Bade if and only if $\liminf_{n\to\infty} \alpha_n < \infty$.
\smallskip
\item{(d)} The algebra $A$ has bounded relative units in the sense of Dales
if and only if the sequence $(\alpha_n)$ is bounded.
}
\medskip
{PROOF:} We know that $A_\alpha$ is a Ditkin algebra (and hence strongly 
regular). Also, for all $x\in \N$, $J_x = M_x$, and $M_x$ has an identity.
Thus, by Proposition 1.3, it is clear that (a), (b), (c) are equivalent,
and that for each of these we only need to check the relevant condition at 
the point $\infty$. Suppose first that 
$\liminf_{n\to\infty} \alpha_n < \infty$. 
Then $(\alpha_n)$ has a bounded subsequence, and so it follows from from
Theorem 2.2 (b) that $M_\infty$ has a bounded approximate identity. Thus
$A_\alpha$ is a strong Ditkin algebra, and has bounded relative units in
the sense of Bade.\par
Conversely, suppose that $A_\alpha$ has bounded relative units in the 
sense of Bade. Then there is a norm bounded sequence of functions 
$(f_n) \subseteq J_\infty$ such that, for all $n\in N$, 
$f_n(\{n\}) = \{1\}$. Set $M= \sup_{n} \norm{f_n}$. We show that
$\liminf_{n\to\infty} \alpha_n \leq M$.
Suppose, for contradiction, that $\liminf_{n\to\infty} \alpha_n > M$. Choose 
$N \in \N$ such that $\inf\{\alpha_j: j \geq N\} > M$.
Then we have
$$\eqalign{1 &= |f_N(N)|\cr
&= \left |{\sum_{j=N}^{\infty}{(f_N(j) - f_N(j+1))}}\right |\cr
&\leq \sum_{j=N}^{\infty}{|(f_N(j) - f_N(j+1))|}\cr
&\leq {{\norm{f_N}}\over{\inf\{\alpha_j: j \geq N\}}} < 1.}$$
This contradiction proves that, as claimed, 
$\liminf_{n\to\infty} \alpha_n \leq M < \infty$. Parts (a) to (c) 
now follow. 
\smallskip
To see (d), suppose first that $(\alpha_n)$ is bounded. 
Set $M=\sup_n \alpha_n$. 
For $n \in \N$, the norm of the identity of $M_n$ is at most $2M+1$.
Also, with $e_n$ as above, $\norm{e_n} \leq M+1$, and so $A_\alpha$
has bounded relative units in the sense of Dales, with bound $2M+1$.
Conversely, since the norm of the identity in $M_n$ is at least 
$\alpha_n$, it is immediate that if $A_\alpha$ has bounded relative units 
in the sense of Dales then $(\alpha_n)$ must be bounded. The result follows.
\qed
\medskip
This result can be used 
to give many elementary examples of strong Ditkin algebras which do not 
have bounded relative units in the sense of Dales.
\medskip
{\bf COROLLARY 2.7.} {\sl Let $\alpha$ be any sequence of positive real
numbers such that $\alpha$ is unbounded but does not diverge to $\infty$.
Then $A_\alpha$ is a strong Ditkin algebra which does not have
bounded relative units in the sense of Dales.}
\smallskip
{PROOF:} This is immediate from Theorem 2.6. \qed
\vfil
\eject
\centerline{\bf REFERENCES}
\medskip
\item{[{\bf 1}]}G.F.~Bachelis, W.A.~Parker, K.A.~Ross, \lq Local units in
$L^1(G)$', {\it Proc. Amer. Math. Soc.} {\bf  31} (1972), 312--313.
\item{[{\bf 2}]}W.G.~Bade, \lq Open problems', in {\it Conference on
Automatic Continuity and Banach Algebras}, 
Proceedings of the Centre for Mathematical Analysis {\bf 21}
(Australian National University, 1989).
\item{[{\bf 3}]}W.G.~Bade, H.G.~Dales, 
\lq The Wedderburn decomposability of some commutative Banach algebras', 
{\it J. Funct. Anal.} {\bf 107} (1992), no. 1, 105--121.
\item{[{\bf 4}]}J.F.~Feinstein, \lq A note on strong Ditkin algebras',
{\it Bull. Austral. Math. Soc.\/} {\bf 52} (1995), 25--30.
\item{[{\bf 5}]}J.F.~Feinstein, \lq Regularity conditions for Banach function 
algebras', Function spaces (Edwardsville, Il, 1994), 
{\it Lecture Notes in Pure and Appl. Math.}, vol. 172, Dekker, New York 1995,
pages 117--122.
\item{[{\bf 6}]}M.M.~Neumann, \lq Commutative Banach algebras and decomposable 
operators', {\it Monatshefte Math.}, 113 (1992) 227-243.
\medskip
\centerline{1991 Mathematics Subject Classification 46J10}
\end